\documentclass[11pt]{article}
\usepackage{amsmath, amssymb, amsthm}
\usepackage{fullpage}
\usepackage{hyperref}
\usepackage{enumitem}
\usepackage{graphicx}
\usepackage{authblk}
\usepackage{pgfplots}
\usepgfplotslibrary{fillbetween}
\pgfplotsset{compat=1.16}

\newtheorem{theorem}{Theorem}[section]
\newtheorem{corollary}[theorem]{Corollary}
\newtheorem{lemma}[theorem]{Lemma}
\newtheorem{observation}[theorem]{Observation}

\theoremstyle{definition}

\newtheorem{construction}[theorem]{Construction}
\theoremstyle{remark}
\newtheorem{remark}[theorem]{Remark}

\newcommand{\BDR}{\textnormal{\textsc{BDR}}}
\newcommand{\BDRab}{\textnormal{\textsc{BDR}}_{\alpha,\beta}}
\newcommand{\Dset}{\mathcal{D}}

\title{On the Complexity of Bipartite Degree Realizability}
\author[$1,2$]{Istv\'an Mikl\'os}

\affil[$1$]{HUN-REN Alfr\'ed R\'enyi Institute of
Mathematics, Budapest, Hungary} 
\affil[$2$]{HUN-REN SZTAKI\\ 1111 Budapest, L\'agym\'anyosi u. 11\\ Hungary}

\begin{document}

\maketitle

\begin{abstract}
We study the \emph{Bipartite Degree Realization} (\BDR) problem: given a degree sequence $D$, decide whether it admits a realization as a bipartite graph. While bipartite realizability for a fixed vertex partition can be decided in polynomial time via the Gale--Ryser theorem, the computational complexity of \BDR{} without a prescribed partition remains unresolved. We address this question through a parameterized analysis.

For constants $0 \le \alpha \le \beta \le 1$, we define $\BDRab$ as the restriction of \BDR{} to degree sequences of length $n$ whose degrees lie in the interval $[\alpha n, \beta n]$. Our main result shows that whenever $\alpha(1-\beta)\ge \beta^2$, or equivalently, $\beta \le \frac{\sqrt{\alpha(\alpha+4)}-\alpha}{2}$, \emph{every} partition of such a sequence into two parts of equal sum is bigraphic. Consequently, $\BDRab$ is solvable in polynomial time for these parameters, by a bounded subset-sum computation. The condition is best possible for this structural statement. The proof is a direct estimate of the Gale--Ryser inequalities.

We further show that, assuming the NP-completeness of unrestricted \BDR, the problem $\BDRab$ remains NP-complete for all $0 < \beta < \frac{1}{2}$ and $0\le \alpha < 1 - \beta - \sqrt{1-2\beta}$, that is, whenever $\sqrt{2\alpha}+\sqrt{1-2\beta}<1$. The proof is a polynomial-time padding reduction. Our results clarify the algorithmic landscape of bipartite degree realization and contribute to the broader study of potentially bipartite degree sequences.
\end{abstract}

\section{Introduction}

A nonnegative integer sequence $D=(d_1,\dots,d_n)$ is called \emph{graphic} if there exists a simple graph on $n$ vertices whose degrees are exactly $d_1,\dots,d_n$. We call such a graph a \emph{realization} of $D$. Graphic sequences were characterized by Erd{\H o}s and Gallai \cite{eg}. Similarly, a pair of degree sequences $(D_1,D_2)$ is \emph{bigraphic} if it has a bipartite graph realization $G=(U,V,E)$ such that the degrees in $U$ are $D_1$ and the degrees in $V$ are $D_2$.
In 1981, Rao surveyed the theory of potentially and forcibly $\mathcal{P}$-graphic degree sequences \cite{rao}. Given a graph property $\mathcal{P}$ and a degree sequence $D$, we say that $D$ is potentially $\mathcal{P}$-graphic if $D$ has a simple graph realization with property $\mathcal{P}$. Similarly, $D$ is forcibly $\mathcal{P}$-graphic if every realization of $D$ has property $\mathcal{P}$. For example, if the sum of the degrees in a degree sequence $D$ of length $n$ is less than $2n-2$, then $D$ is forcibly disconnected. If the sum is at least $2n-2$, and there is no $0$ degree, then $D$ is potentially connected. If the minimum degree is greater than $\left\lfloor\frac{n-2}{2}\right\rfloor$, then $D$ is forcibly connected.

One of the open problems listed by Rao concerns potentially bipartite degree sequences. 
Following his terminology, we call a degree sequence \emph{potentially bipartite} if it admits at least one bipartite realization, and \emph{forcibly not bipartite} if none of its realizations is bipartite. The complementary notion of \emph{forcibly} bipartite sequences was recently characterized completely by Bar-Noy, B\"ohnlein, Peleg and Rawitz \cite{bbpr-forcibly}.
We will denote the potentially bipartite degree sequence problem as the \BDR{} problem, that is, the Bipartite Degree Realization problem. The input of the \BDR{} problem is a degree sequence $D$, and the problem is to decide whether $D$ has a bipartite graph realization. The computational complexity of this problem is still unknown, although the following facts are well known. The Gale--Ryser theorem provides a polynomial-time algorithm to decide whether a pair of degree sequences $(D_1,D_2)$ has a bipartite graph realization $G = (U,V,E)$ such that the degrees in $U$ are $D_1$ and the degrees in $V$ are $D_2$ \cite{gale, ryser}. Also, there is a polynomial-time dynamic programming algorithm to decide whether a degree sequence $D$ can be split into two degree sequences $D_1$ and $D_2$ with equal sums, which is a necessary condition for having a bipartite graph realization (the degrees of a graphic sequence of length $n$ are smaller than $n$, so this subset-sum instance has polynomially bounded numbers). However, there might be an exponential number of ways to split $D$ into two degree sequences with equal sums \cite{barnoy1}, and it is not clear how to decide whether any of them satisfies the necessary and sufficient conditions given by Gale and Ryser. 

In this paper, we consider a parametric version of \BDR. The inputs of $\BDRab$ are degree sequences of length $n$ in which every degree is between $\alpha n$ and $\beta n$. The problem is still to decide whether a bipartite realization exists. We obtain the following results.
\begin{itemize}
    \item If $\alpha(1-\beta)\ge\beta^2$, equivalently, $\beta\le \frac{\sqrt{\alpha(\alpha+4)}-\alpha}{2}$, then every split of a degree sequence of length $n$ with degrees in $[\alpha n,\beta n]$ into two subsequences with equal sums is bigraphic (Theorem~\ref{thm:main-structural}). In particular, such a sequence is potentially bipartite if and only if it can be split into two subsequences with equal sums, and thus $\BDRab$ is in P for these parameters (Corollary~\ref{cor:poly}). The condition $\alpha(1-\beta)\ge\beta^2$ cannot be weakened in this structural statement (Remark~\ref{rem:sharp}).
    \item If \BDR{} is NP-complete, then $\BDRab$ is NP-complete for all $0<\beta<\frac12$ and $0\le\alpha<1-\beta-\sqrt{1-2\beta}$ (Theorem~\ref{thm:hard}). We will show that, for $0\le \alpha<\frac12$ and $0<\beta<\frac12$, this condition is equivalent to $\beta>\sqrt{2\alpha}-\alpha$, and also to $\sqrt{2\alpha}+\sqrt{1-2\beta}<1$ (Lemma~\ref{lem:hard-alg}). This result follows from a polynomial-time padding reduction that preserves bipartite realizability.
\end{itemize}
The two parameter regions are depicted in Figure~\ref{fig:regions}. We remark that the case $\alpha>\frac12$ is trivial and is not part of our results: if every degree is larger than $\frac{n}{2}$, then $D$ is forcibly not bipartite, since in any bipartite graph on $n$ vertices one of the two vertex classes has at most $\frac{n}{2}$ vertices, and the degrees in the other vertex class are bounded by its size.

\begin{figure}[t]
\centering
\begin{tikzpicture}
\begin{axis}[
    width=0.62\textwidth, height=0.55\textwidth,
    xmin=0, xmax=0.5, ymin=0, ymax=0.5,
    xlabel={$\alpha$}, ylabel={$\beta$},
    xtick={0,0.1,0.2,0.3,0.4,0.5}, ytick={0,0.1,0.2,0.3,0.4,0.5},
    tick label style={font=\small},
    legend style={font=\small, at={(0.98,0.03)}, anchor=south east},
    legend cell align=left,
]
\addplot[name path=diag, domain=0:0.5, samples=2, black, thin, dashed] {x};
\addlegendentry{$\beta=\alpha$}
\addplot[name path=T, domain=0:0.5, samples=400, black, thick] {(sqrt(x*(x+4))-x)/2};
\addlegendentry{$\beta=\frac{\sqrt{\alpha(\alpha+4)}-\alpha}{2}$}
\addplot[name path=H, domain=0:0.5, samples=400, black, thick, densely dotted] {sqrt(2*x)-x};
\addlegendentry{$\beta=\sqrt{2\alpha}-\alpha$}
\addplot[name path=top, domain=0:0.5, samples=2, draw=none, forget plot] {0.5};
\addplot[gray!45, forget plot] fill between[of=T and diag];
\addplot[gray!15, forget plot] fill between[of=top and H];
\node[font=\small] at (axis cs:0.30,0.36) {P};
\node[font=\small, align=center] at (axis cs:0.075,0.44) {NP-c.\ if \BDR{} is};
\node[font=\small] at (axis cs:0.13,0.34) {?};
\end{axis}
\end{tikzpicture}
\caption{The parameter plane of $\BDRab$ for $0\le\alpha\le\beta\le\frac12$. In the dark region ($\alpha\le\beta\le\frac{\sqrt{\alpha(\alpha+4)}-\alpha}{2}$) every equal-sum split is bigraphic and the problem is in P (Theorem~\ref{thm:main-structural}). In the light region ($\sqrt{2\alpha}-\alpha<\beta<\frac12$) the problem is NP-complete provided \BDR{} is NP-complete (Theorem~\ref{thm:hard}). The complexity between the two curves is open.}
\label{fig:regions}
\end{figure}

The computational complexity of deciding whether a degree sequence admits a bipartite realization has recently been studied intensively by Bar-Noy, B\"ohnlein, Peleg and Rawitz, both for simple graphs and for multigraphs \cite{barnoy1, barnoy2, bbpr-equal, bbpr-forcibly}.
Our work is concerned exclusively with simple graphs. 
Accordingly, we restrict attention to those of their results that apply to simple bipartite graphs and are directly comparable to ours.
In this restricted setting, Bar-Noy et al.\ identify polynomial-time solvable cases primarily under strong structural assumptions on the degree sequence, most notably when the sequence has a suitable high-low partition \cite{barnoy1}. A high-low partition of a degree sequence $D$ is a partition into $D_1$ and $D_2$ such that for all $d_i\in D_1$ and $d_j\in D_2$, $d_i\le d_j$. A high-low partition $(D_1,D_2)$ is \emph{well-behaved} if the maximum degree in each vertex class is at most the size of the other vertex class. One of the main results of \cite{barnoy1} is the following. If $D$ has a well-behaved high-low partition, then $D$ has a bipartite graph realization if and only if the well-behaved high-low partition has a bipartite graph realization. That is, if the well-behaved high-low partition does not have a bipartite graph realization, then none of the equal-sum partitions of the underlying degree sequence is bigraphic.
These results assume a special global property (namely, a well-behaved high-low partition) rather than extremal degree constraints. In a related work \cite{bbpr-equal}, the same authors study \emph{even} sequences, in which every value occurs with even multiplicity, so that the sequence can be split into two identical halves (the \emph{equal partition}); they prove that a graphic even sequence is either bigraphic, or it can be realized by a bipartite multigraph with edge multiplicities at most $2$.

In contrast, in the present work, we identify degree sequence classes in which \emph{every} partition into two parts with equal sums has a bipartite graph realization. This establishes polynomial-time decidability for bipartite realizability under explicit extremal degree conditions of the form \(\alpha n \le d_i \le \beta n\), covering a broad class of dense degree sequences with linear minimum and maximum degree bounds. Our results apply also when a degree sequence does not have a high-low partition of the above kind, and also
when the degree sequence exhibits large variance, provided all degrees remain within the prescribed interval. 
The largest difference between the maximum and minimum degrees that our tractability result allows is obtained by maximizing $\frac{\sqrt{\alpha(\alpha+4)}-\alpha}{2}n -\alpha n$, which turns out to be $(3-2\sqrt{2})n\approx 0.17157n$, attained at $\alpha = \frac{3\sqrt{2}}{2}-2$.
This places our tractability region in a different part of the parameter space than the one considered by Bar-Noy et al.\ for simple graphs.

In another work, Bar-Noy et al.\ proved that the corresponding bipartite realization problem for multigraphs is NP-hard \cite{barnoy2}. In that setting, degrees are allowed to be arbitrarily large and are not bounded as a function of the number of vertices. From this perspective, the hardness result aligns with the inherent difficulty of partitioning sequences of large integers into two parts with equal sums, a problem that is NP-complete \cite{garey-johnson}. Consequently, these multigraph results are not directly comparable to the simple graph setting considered here.

The computational complexity of \BDR{} for simple graphs, where degrees are bounded by $n-1$, remains open \cite{barnoy1,barnoy2}. Our results contribute to this open direction by showing that if \BDR{} for simple graphs were NP-complete, then NP-completeness would persist even under much stronger restrictions; in particular, it would persist for sequences whose maximum degree is at most $cn$, for any fixed constant $0<c<\frac{1}{2}$ (take $\alpha=0$ and $\beta=c$ in Theorem~\ref{thm:hard}).

The paper is organized as follows. Section~\ref{sec:prelim} collects the definitions and the tools we use. Section~\ref{sec:poly} contains the structural theorem and the polynomial algorithm, and Section~\ref{sec:hard} contains the conditional hardness result. We conclude with some open questions in Section~\ref{sec:concl}.

\section{Preliminaries}\label{sec:prelim}

For constants $0\le\alpha\le\beta\le 1$ and a positive integer $n$, let $\Dset_n(\alpha,\beta)$ denote the set of sequences $D=(d_1,\ldots,d_n)$ of nonnegative integers with $\alpha n\le d_i\le\beta n$ for all $i$. The instances of $\BDRab$ are the sequences in $\bigcup_{n\ge 1}\Dset_n(\alpha,\beta)$, and the question is whether the given sequence has a bipartite realization. We stress that $\alpha n$ and $\beta n$ are not required to be integers; they are used only as real lower and upper bounds for the (integer) degrees. 

A \emph{split} of a sequence $D$ is a partition of its entries into two subsequences $D_1$ and $D_2$, and the split is an \emph{equal-sum split} if the entries of $D_1$ and the entries of $D_2$ have the same sum. Clearly, $D$ is potentially bipartite if and only if it has a split $(D_1,D_2)$ that is bigraphic, and every bigraphic split is an equal-sum split.

The bigraphicality of a pair of sequences can be checked with the Gale--Ryser theorem \cite{gale,ryser}.
\begin{theorem}[Gale--Ryser]\label{thm:gale-ryser}
Let $D_1 = (d_{1,1}, d_{1,2}, \ldots, d_{1,s_1})$ and $D_2 = (d_{2,1},d_{2,2},\ldots,d_{2,s_2})$ be sequences of nonnegative integers, where $d_{1,1} \ge d_{1,2}\ge \cdots \ge d_{1,s_1}$. Then $(D_1,D_2)$ is bigraphic if and only if $\sum_{i=1}^{s_1} d_{1,i}=\sum_{j=1}^{s_2} d_{2,j}$ and for all $k=1,2,\ldots,s_1$,
\begin{equation}
\sum_{i=1}^k d_{1,i} \le \sum_{j=1}^{s_2} \min\{k,d_{2,j}\}.\label{eq:gale-ryser}    
\end{equation}
\end{theorem}
Note that only the entries of $D_1$ have to be ordered: the right-hand side of \eqref{eq:gale-ryser} does not depend on the order of the entries of $D_2$. Note also that the roles of $D_1$ and $D_2$ may be interchanged; it is enough to verify \eqref{eq:gale-ryser} for one of the two possible choices.

The necessity of \eqref{eq:gale-ryser} is a simple counting argument, which we will also need in a slightly more general form, together with the case of equality.
\begin{observation}\label{obs:counting}
Let $G=(X,Y,E)$ be a bipartite graph and let $A\subseteq X$. Then
\begin{equation}
\sum_{v\in A}\deg(v) \le \sum_{y\in Y}\min\{|A|,\deg(y)\}.\label{eq:counting}
\end{equation}
Moreover, equality holds if and only if every $y\in Y$ with $\deg(y)\ge |A|$ is adjacent to all vertices of $A$, and every $y\in Y$ with $\deg(y)<|A|$ has all its neighbors in $A$.
\end{observation}
\begin{proof}
Since every edge incident with a vertex of $A$ has its other endpoint in $Y$, the left-hand side of \eqref{eq:counting} equals $\sum_{y\in Y}|N(y)\cap A|$. For every $y\in Y$, we have $|N(y)\cap A|\le\min\{|A|,\deg(y)\}$, and equality holds if and only if $N(y)\supseteq A$ (when $\deg(y)\ge|A|$) or $N(y)\subseteq A$ (when $\deg(y)<|A|$). Summing over $y\in Y$ gives both statements.
\end{proof}

We will also need the following folklore fact on semi-regular bipartite graphs. We include a proof for completeness.
\begin{lemma}\label{lem:biregular}
Let $x,y,r,c$ be positive integers with $xr=yc$, $r\le y$ and $c\le x$. Then there exists a simple bipartite graph with vertex classes $X$ and $Y$, $|X|=x$ and $|Y|=y$, in which every vertex of $X$ has degree $r$ and every vertex of $Y$ has degree $c$.
\end{lemma}
\begin{proof}
Let $X=\{0,1,\ldots,x-1\}$ and $Y=\mathbb{Z}_y$, and connect $i\in X$ to the vertices $ir, ir+1,\ldots, ir+r-1$ (modulo $y$). Since $r\le y$, these are $r$ distinct vertices, so the graph is simple and every vertex of $X$ has degree $r$. The list of all edge endpoints in $Y$, taken in the order $i=0,1,\ldots,x-1$, is the sequence $0,1,\ldots,xr-1$ reduced modulo $y$. Since $xr=yc$, every residue occurs exactly $c$ times, so every vertex of $Y$ has degree $c$.
\end{proof}

\section{\texorpdfstring{Polynomial algorithm for $\BDRab$}{Polynomial algorithm for BDR(alpha,beta)}}\label{sec:poly}

We first state the equivalent forms of the condition of our tractability result.
\begin{lemma}\label{lem:tractable-alg}
Let $\alpha \ge 0$ and $0 \le \beta < 1$. Then the following are equivalent:
\begin{enumerate}[label=(\roman*)]
    \item $\alpha \ge \frac{\beta^2}{1-\beta}$;
    \item $\alpha(1-\beta)\ge\beta^2$;
    \item $\beta \le \frac{\sqrt{\alpha(\alpha+4)} - \alpha}{2}$.
\end{enumerate}
Moreover, if $\alpha+\beta>0$, then these conditions are also equivalent to
\begin{enumerate}[label=(\roman*)]\setcounter{enumi}{3}
    \item $\beta\le\frac{\alpha}{\alpha+\beta}$.
\end{enumerate}
Finally, if $0\le \alpha\le\beta$ and these conditions hold, then $\beta\le\frac12$.
\end{lemma}
\begin{proof}
Since $1-\beta>0$, (i) and (ii) are equivalent. Condition (ii) can be rearranged as
$$
\beta^2+\alpha\beta-\alpha\le 0 .
$$
The quadratic polynomial $t^2+\alpha t-\alpha$ has the two real roots
$$
t_{\pm}=\frac{-\alpha\pm\sqrt{\alpha(\alpha+4)}}{2},
$$
and it is nonpositive exactly on the interval $[t_-,t_+]$. Since $t_-\le 0\le\beta$, the inequality $\beta^2+\alpha\beta-\alpha\le0$ holds if and only if $\beta\le t_+$, which is (iii). If $\alpha+\beta>0$, then multiplying (iv) by $\alpha+\beta$ gives $\alpha\beta+\beta^2\le\alpha$, which is again (ii). Finally, if $\alpha\le\beta$ and (ii) holds, then $\beta^2\le\alpha(1-\beta)\le\beta(1-\beta)$, hence $\beta=0$ or $\beta\le 1-\beta$; in both cases $\beta\le\frac12$.
\end{proof}

The main result of this section is the following structural theorem.
\begin{theorem}\label{thm:main-structural}
Let $0\le\alpha\le\beta$ satisfy $\alpha(1-\beta)\ge\beta^2$, that is, $\beta\le\frac{\sqrt{\alpha(\alpha+4)}-\alpha}{2}$. Let $D\in\Dset_n(\alpha,\beta)$. Then every equal-sum split $(D_1,D_2)$ of $D$ is bigraphic.
\end{theorem}
\begin{proof}
If $\alpha=0$, then the condition $\alpha(1-\beta)\ge\beta^2$ forces $\beta=0$, so all degrees are $0$ and the statement is trivial. Hence we may assume $\alpha>0$, and then $0<\alpha\le\beta\le\frac12$ by Lemma~\ref{lem:tractable-alg}.

Let $(D_1,D_2)$ be an equal-sum split of $D$, let $s_1$ and $s_2$ be the lengths of $D_1$ and $D_2$, respectively, and let $\Sigma$ be their common sum. We have $s_1+s_2=n$, and we may assume without loss of generality that
$$
s_1\le s_2 .
$$
We write the entries of $D_1$ in non-increasing order as $d_{1,1}\ge\cdots\ge d_{1,s_1}$, and the entries of $D_2$ in any order as $d_{2,1},\ldots,d_{2,s_2}$. Let
$$
a:=\alpha n,\qquad b:=\beta n,\qquad \delta:=b-a\ge0 .
$$
These are real numbers, and every entry of $D_1$ and $D_2$ lies in the interval $[a,b]$.

\medskip\noindent\textbf{Step 1: both classes are large.} Since every entry of $D_1$ is at most $b$ and every entry of $D_2$ is at least $a$, we have $s_2a\le\Sigma\le s_1b$; hence $(n-s_1)a\le s_1b$, that is, $s_1\ge\frac{a}{a+b}n$. By the same argument, $s_2\ge\frac{a}{a+b}n$. By Lemma~\ref{lem:tractable-alg}(iv),
\begin{equation}
\min\{s_1,s_2\}\;\ge\;\frac{\alpha}{\alpha+\beta}\,n\;\ge\;\beta n=b.\label{eq:large-classes}
\end{equation}

\medskip\noindent\textbf{Step 2: the Gale--Ryser inequalities.}
By Theorem~\ref{thm:gale-ryser}, it suffices to show that
$$
\mathrm{L}(k):=\sum_{i=1}^k d_{1,i}\;\le\;\mathrm{R}(k):=\sum_{j=1}^{s_2}\min\{k,d_{2,j}\}
$$
for every integer $1\le k\le s_1$. We distinguish three cases.

\smallskip\noindent\emph{Case (i): $k\le a$.} Then $d_{2,j}\ge a\ge k$ for every $j$, hence $\mathrm{R}(k)=s_2k$. On the other hand, $\mathrm{L}(k)\le kb\le ks_2$ by \eqref{eq:large-classes}.

\smallskip\noindent\emph{Case (ii): $k\ge b$.} Then $d_{2,j}\le b\le k$ for every $j$, hence $\mathrm{R}(k)=\Sigma\ge \mathrm{L}(k)$.

\smallskip\noindent\emph{Case (iii): $a<k<b$.} In particular, $\delta>0$. We bound the two sides separately. Since the largest $k$ entries of $D_1$ are at most $b$ and the remaining $s_1-k$ entries are at least $a$, we have
\begin{equation}
\mathrm{L}(k)\le\min\{kb,\;\Sigma-(s_1-k)a\}.\label{eq:L-bound}
\end{equation}
For the right-hand side, we use that for every real $d\in[a,b]$,
$$
\min\{k,d\}\;\ge\;a+(d-a)\frac{k-a}{b-a}.
$$
Indeed, if $d\le k$, this is equivalent to $(d-a)\bigl(1-\frac{k-a}{b-a}\bigr)\ge0$, and if $d>k$, it is equivalent to $(k-a)\bigl(1-\frac{d-a}{b-a}\bigr)\ge 0$; both hold since $a\le d\le b$ and $a<k<b$. (Geometrically, the concave function $d\mapsto\min\{k,d\}$ lies above its chord over $[a,b]$.) Summing over the entries of $D_2$, we obtain
\begin{equation}
\mathrm{R}(k)\;\ge\;s_2a+(\Sigma-s_2a)\frac{k-a}{b-a}.\label{eq:R-bound}
\end{equation}
For real $x\in[a,b]$, define
$$
\Phi(x):=s_2a+(\Sigma-s_2a)\frac{x-a}{b-a}-\min\{xb,\;\Sigma-(s_1-x)a\}.
$$
By \eqref{eq:L-bound} and \eqref{eq:R-bound}, $\mathrm{R}(k)-\mathrm{L}(k)\ge\Phi(k)$, so it is enough to show that $\Phi(x)\ge0$ for all $x\in[a,b]$. Although only the values of $\Phi$ at integers $k$ with $a<k<b$ are needed, we consider $\Phi$ as a function of a real variable on the closed interval $[a,b]$: the minimum of a convex piecewise linear function on a closed interval is attained at one of the endpoints or at a breakpoint, and this makes it easy to locate the minimum. Note that $\Phi(a)$ and $\Phi(b)$ are used only as values of this auxiliary function; they need not be related to any Gale--Ryser inequality, and $a$ and $b$ need not be integers. The function $\Phi$ is the difference of an affine function and the minimum of two affine functions; hence it is convex and piecewise linear, with at most one breakpoint, namely at the point $x^*$ where $x^*b=\Sigma-(s_1-x^*)a$, that is,
$$
x^*=\frac{\Sigma-s_1a}{b-a}.
$$
Therefore, the minimum of $\Phi$ on $[a,b]$ is attained at $x=a$, at $x=b$, or at $x=x^*$ if $a<x^*<b$. We check these points.

At $x=a$, using \eqref{eq:large-classes},
$$
\Phi(a)\ge s_2a-ab=a(s_2-b)\ge0 .
$$
At $x=b$, using \eqref{eq:large-classes} again,
$$
\Phi(b)\ge\Sigma-\bigl(\Sigma-(s_1-b)a\bigr)=a(s_1-b)\ge0 .
$$
Now assume that $a<x^*<b$, and write $x^*=a+\tau\delta$ with $0<\tau<1$. From the definition of $x^*$, we have $\Sigma=s_1a+x^*\delta=s_1a+a\delta+\tau\delta^2$, and $\min\{x^*b,\Sigma-(s_1-x^*)a\}=x^*b$. Therefore,
\begin{align}
\Phi(x^*)&=s_2a+\tau(\Sigma-s_2a)-(a+\tau\delta)b\nonumber\\
&=(1-\tau)s_2a+\tau(s_1a+a\delta+\tau\delta^2)-ab-\tau\delta b\nonumber\\
&=a\bigl[(1-\tau)s_2+\tau s_1-b\bigr]-\tau(1-\tau)\delta^2,\label{eq:phi-star}
\end{align}
where in the last step we used $\tau a\delta-\tau\delta b=-\tau\delta^2$. We consider two subcases.

\smallskip\noindent\emph{Subcase $\tau\le\frac12$.} Since $s_1\le s_2$,
$$
(1-\tau)s_2+\tau s_1=\frac{s_1+s_2}{2}+\Bigl(\frac12-\tau\Bigr)(s_2-s_1)\ge\frac{n}{2}.
$$
Using also $\tau(1-\tau)\le\frac14$, \eqref{eq:phi-star} gives
$$
\Phi(x^*)\ge a\Bigl(\frac n2-b\Bigr)-\frac{\delta^2}{4}
=\frac{n^2}{4}\Bigl(2\alpha(1-\beta)-\alpha^2-\beta^2\Bigr)
\ge\frac{n^2}{4}\Bigl(2\beta^2-\alpha^2-\beta^2\Bigr)=\frac{n^2}{4}(\beta^2-\alpha^2)\ge0,
$$
where we used $\alpha(1-\beta)\ge\beta^2$ and $\alpha\le\beta$.

\smallskip\noindent\emph{Subcase $\tau>\frac12$.} Here we use the lower bound $\Sigma\ge s_2a=(n-s_1)a$ on the sum of $D_2$. Together with $\Sigma=s_1a+a\delta+\tau\delta^2$, it gives
$$
2s_1a\ge na-a\delta-\tau\delta^2.
$$
Writing $(1-\tau)s_2+\tau s_1=(1-\tau)n+(2\tau-1)s_1$ and using $2\tau-1>0$, we get from \eqref{eq:phi-star}
\begin{align*}
\Phi(x^*)&=a(1-\tau)n+(2\tau-1)as_1-ab-\tau(1-\tau)\delta^2\\
&\ge a(1-\tau)n+\frac{2\tau-1}{2}\bigl(na-a\delta-\tau\delta^2\bigr)-ab-\tau(1-\tau)\delta^2\\
&=\frac{an}{2}-ab-\frac{(2\tau-1)a\delta}{2}-\frac{\tau\delta^2}{2}=:\psi(\tau).
\end{align*}
The function $\psi$ is affine and non-increasing in $\tau$, hence for $\frac12<\tau<1$,
$$
\Phi(x^*)\ge\psi(1)=\frac{an}{2}-ab-\frac{a\delta}{2}-\frac{\delta^2}{2}=\frac{an-ab-b^2}{2}=\frac{n^2}{2}\Bigl(\alpha(1-\beta)-\beta^2\Bigr)\ge0 .
$$

We conclude that $\Phi\ge0$ on $[a,b]$; hence $\mathrm{L}(k)\le\mathrm{R}(k)$ also in Case (iii). Thus all Gale--Ryser inequalities hold, and $(D_1,D_2)$ is bigraphic.
\end{proof}

\begin{remark}
The bounds \eqref{eq:L-bound} and \eqref{eq:R-bound} are attained when $D_1$ and $D_2$ are \emph{least balanced} sequences, that is, when all their entries are equal to $a$ or $b$, except at most one entry. This is in accordance with the well-known fact that, for fixed length, sum, minimum and maximum, the least balanced degree sequences are the hardest to realize (see, e.g., \cite{at94, limiklos}). The proof above shows this directly, without referring to these results.
\end{remark}

\begin{corollary}\label{cor:poly}
Let $0\le\alpha\le\beta$ satisfy $\alpha(1-\beta)\ge\beta^2$. Then a sequence $D\in\Dset_n(\alpha,\beta)$ is potentially bipartite if and only if it has an equal-sum split. Consequently, $\BDRab$ can be solved in $O(n^3)$ time, and if a bipartite realization exists, then one can be constructed in polynomial time.
\end{corollary}
\begin{proof}
Every bigraphic split is an equal-sum split, and conversely, every equal-sum split is bigraphic by Theorem~\ref{thm:main-structural}. Whether $D$ has an equal-sum split is a subset-sum question: the sum $S$ of the degrees must be even, and some sub-multiset of the degrees must have sum $\frac{S}{2}$. Since $S\le \beta n^2\le n^2$, the standard dynamic programming algorithm for subset sum decides this in $O(nS)=O(n^3)$ time and also returns an equal-sum split if there is one. Given a bigraphic pair, a realization can be constructed in polynomial time, for example, by the constructive proof of the Gale--Ryser theorem \cite{ryser}.
\end{proof}

\begin{remark}\label{rem:sharp}
The condition $\alpha(1-\beta)\ge\beta^2$ cannot be weakened in Theorem~\ref{thm:main-structural}. Assume that $0<\alpha<\beta$ and $\alpha(1-\beta)<\beta^2$, that is, $\frac{\alpha}{\alpha+\beta}<\beta$. Let $a_0:=\lceil\alpha n\rceil$, $b_0:=\lfloor\beta n\rfloor$ and $p:=\left\lceil\frac{(n-1)a_0+b_0}{a_0+b_0}\right\rceil$. Since $a_0=\alpha n+O(1)$ and $b_0=\beta n+O(1)$, we have $p=\frac{\alpha}{\alpha+\beta}n+O(1)$, and hence, for all sufficiently large $n$, $a_0\le b_0$ and $p<b_0$ (and also $p<\frac n2$). Consider a split $(D_1,D_2)$ in which $D_1$ has $p$ entries from $[a_0,b_0]$, and $D_2$ has $n-p$ entries, one of them equal to $b_0$ and the others from $[a_0,b_0]$. The possible sums of $D_1$ form the integer interval $[pa_0,pb_0]$, and the possible sums of $D_2$ form the integer interval $[b_0+(n-p-1)a_0,(n-p)b_0]$. These two intervals intersect, since $pa_0\le(n-p)b_0$, and $b_0+(n-p-1)a_0\le pb_0$ by the choice of $p$; so the entries can be chosen to obtain an equal-sum split. The resulting sequence is in $\Dset_n(\alpha,\beta)$, but the split is not bigraphic, since $D_2$ contains a degree $b_0$, which is larger than the number $p$ of entries of $D_1$.

On the other hand, this does not mean that $\BDRab$ is hard for such parameters: it only means that, beyond this threshold, the existence of an equal-sum split is no longer sufficient for potential bipartiteness.
\end{remark}

\section{The conditional NP-completeness result}\label{sec:hard}

Now we turn to the other main result of the paper. 
We start with the equivalent forms of the condition of the hardness result. The key observation is the identity
\begin{equation}
1-\beta-\sqrt{1-2\beta}=\frac{\left(1-\sqrt{1-2\beta}\right)^2}{2},\label{eq:identity}
\end{equation}
which follows from $1-\beta=\frac{1+(1-2\beta)}{2}$ by expanding the square on the right-hand side.
\begin{lemma}\label{lem:hard-alg}
Let $0\le\alpha$ and $0<\beta<\frac12$. Then the following are equivalent:
\begin{enumerate}[label=(\roman*)]
    \item $\alpha<1-\beta-\sqrt{1-2\beta}$;
    \item $\sqrt{2\alpha}+\sqrt{1-2\beta}<1$;
    \item $\alpha<\frac12$ and $\beta>\sqrt{2\alpha}-\alpha$;
    \item $\alpha<\frac12$ and $(\alpha+\beta)^2>2\alpha$.
\end{enumerate}
Moreover, if these conditions hold, then $\alpha<\beta$.
\end{lemma}
\begin{proof}
Let $u:=\sqrt{2\alpha}\ge 0$ and $v:=\sqrt{1-2\beta}\in(0,1)$; then $\alpha=\frac{u^2}{2}$ and $\beta=\frac{1-v^2}{2}$.

(i)$\Leftrightarrow$(ii): By \eqref{eq:identity}, (i) is equivalent to $u^2<(1-v)^2$, and since $u\ge 0$ and $1-v>0$, this is equivalent to $u<1-v$.

(ii)$\Leftrightarrow$(iii): Both (ii) and (iii) imply $u<1$, that is, $\alpha<\frac12$; so assume $u<1$. Then (iii) reads $\frac{1-v^2}{2}>u-\frac{u^2}{2}$, which is equivalent to $(1-u)^2>v^2$, and since $1-u>0$ and $v>0$, to $1-u>v$.

(iii)$\Leftrightarrow$(iv): Since $\alpha+\beta>0$, we have $(\alpha+\beta)^2>2\alpha$ if and only if $\alpha+\beta>\sqrt{2\alpha}$.

Finally, $1-2\beta<\sqrt{1-2\beta}$ since $0<1-2\beta<1$; hence $1-\beta-\sqrt{1-2\beta}<\beta$, and (i) implies $\alpha<\beta$.
\end{proof}

We do not know whether or not \BDR{} is NP-complete. However, we can prove the following theorem.
\begin{theorem}\label{thm:hard}
Let $0<\beta<\frac{1}{2}$ and $0\le\alpha<1-\beta-\sqrt{1-2\beta}$. If \BDR{} is NP-complete, then $\BDRab$ is also NP-complete.
\end{theorem}

The rest of this section is devoted to the proof of Theorem~\ref{thm:hard}. Clearly, $\BDRab$ is in NP: a split together with a bipartite realization is a certificate that can be checked in polynomial time. We give a polynomial-time many-one reduction from \BDR{} to $\BDRab$.

\subsection{Choice of the constants}
By Lemma~\ref{lem:hard-alg}, $\sqrt{1-2\beta}<1-\sqrt{2\alpha}$. Fix a rational number $\frac{p}{q}$ with positive integers $p<q$ such that
\begin{equation}
\sqrt{1-2\beta}<\frac{p}{q}<1-\sqrt{2\alpha}. \label{eq:choice-r}
\end{equation}
Define
$$
\gamma_1:=(q-p)^2-2\alpha q^2\qquad\text{and}\qquad\gamma_2:=p^2-(1-2\beta)q^2 .
$$
By \eqref{eq:choice-r}, $\gamma_1>0$ and $\gamma_2>0$. Finally, let $C$ be a positive integer with
\begin{equation}
2C\ge\max\left\{\frac{2\alpha}{\gamma_1},\frac{1-2\beta}{\gamma_2}\right\}.\label{eq:choice-C}
\end{equation}
The integers $p$, $q$ and $C$ depend only on $\alpha$ and $\beta$, and not on the input of the reduction. Theorem~\ref{thm:hard} is a statement about one fixed pair $(\alpha,\beta)$, so it suffices to show that, for this pair, a polynomial-time reduction exists. Such integers $p$, $q$ and $C$ exist, and the reduction simply uses them as built-in constants. 

\begin{remark}
With $r:=\frac{p}{q}$, the numbers $\tilde\alpha:=\frac{(1-r)^2}{2}$ and $\tilde\beta:=\frac{1-r^2}{2}$ are rational, they satisfy $\alpha<\tilde\alpha<\tilde\beta<\beta$ by \eqref{eq:choice-r}, and $\tilde\alpha=1-\tilde\beta-\sqrt{1-2\tilde\beta}$ by \eqref{eq:identity}. That is, the point $(\tilde\alpha,\tilde\beta)$ lies on the boundary curve of Theorem~\ref{thm:hard}, and the gadget below lives in $\Dset(\tilde\alpha,\tilde\beta)$ up to a lower-order padding. On this boundary, the Gale--Ryser inequality belonging to the large vertices of the gadget holds with equality, which forces the structure of every bipartite realization.
\end{remark}

\subsection{Preprocessing}
Let $D=(d_1,\ldots,d_m)$ be an instance of \BDR. We first remove all $0$ entries; this does not change whether $D$ has a bipartite realization, since isolated vertices can be put into either vertex class. If no entry remains, then $D$ is a yes-instance. Let $S:=\sum_{i=1}^m d_i$. If $S$ is odd, or $\max_i d_i\ge m$, or $\max_i d_i>\frac S2$, then $D$ is a no-instance: in the first two cases $D$ is not even graphic, and in the third case, a vertex of degree $d_i$ would need $d_i$ neighbors in the other vertex class, whose degree sum is $\frac{S}{2}$. In these trivial cases, the reduction outputs a fixed yes-instance or a fixed no-instance of $\BDRab$, for example, the output of Construction~\ref{constr} below on the sequence $(1,1)$ (a yes-instance), or on the sequence $(2,2,2)$ (a no-instance, since its only realization is a triangle); the correctness of these choices follows from Lemmas~\ref{lem:forward} and \ref{lem:backward}.

Hence, from now on, we assume that
\begin{equation}
1\le d_i\le\frac{S}{2}\ \text{ for all } i,\qquad S \text{ is even},\qquad S\ge 2, \qquad S<m^2.\label{eq:assumptions}
\end{equation}
Note that $m\le S$, since every entry is at least $1$.

\subsection{The construction}
Let
$$
K:=2CS^2,
$$
and define the integers
\begin{align*}
a&:=(q-p)^2K, & L&:=q(q-p)K, & s&:=pqK,\\
c&:=p(q-p)K, & b&:=(q^2-p^2)K, & W&:=b+S.
\end{align*}
These numbers satisfy the identities
\begin{align}
L+c&=b, \label{eq:I1}\\
Lc&=sa, \label{eq:I2}\\
Lb&=L^2+sa, \label{eq:I3}
\end{align}
since $q(q-p)+p(q-p)=(q-p)(q+p)$, $q(q-p)\cdot p(q-p)=pq\cdot(q-p)^2$, and \eqref{eq:I3} follows from \eqref{eq:I1} and \eqref{eq:I2}. Furthermore,
\begin{equation}
a<L\le b<W,\qquad c\le s.\label{eq:ineq-param}
\end{equation}
Indeed, $(q-p)^2<q(q-p)$ since $p>0$, $L\le b$ follows from \eqref{eq:I1}, and $p(q-p)\le pq$.

\begin{construction}\label{constr}
Given $D=(d_1,\ldots,d_m)$ satisfying \eqref{eq:assumptions}, let $D'$ be the following sequence of length
$$
N:=2L+2S+2s=2q^2K+2S .
$$
\begin{itemize}
\item $2L$ entries equal to $W$; the corresponding vertices are called \emph{large vertices};
\item $2S$ entries $L+d'_1, L+d'_2,\ldots,L+d'_{2S}$, where $(d'_1,\ldots,d'_{2S})=(d_1,\ldots,d_m,0,\ldots,0)$ is $D$ padded with $2S-m$ zeros; the corresponding vertices are called \emph{original vertices}, and for an original vertex $v$ with entry $L+d'_i$ we write $\hat d(v):=d'_i$;
\item $2s$ entries equal to $a$; the corresponding vertices are called \emph{small vertices}.
\end{itemize}
\end{construction}

\begin{lemma}\label{lem:size}
The sequence $D'$ can be computed from $D$ in polynomial time, and $D'\in\Dset_N(\alpha,\beta)$.
\end{lemma}
\begin{proof}
Since $p$, $q$, $C$ are constants and $S<m^2$ by \eqref{eq:assumptions}, we have $N=O(S^2)=O(m^4)$, and all entries of $D'$ are at most $N$; hence $D'$ can be computed in polynomial time.

By \eqref{eq:ineq-param} and $\hat d(v)\le\frac{S}{2}$, every entry of $D'$ is at least $a$ and at most $W$ (note that $L+\frac S2<b+S=W$). Hence it is enough to show that $a\ge\alpha N$ and $W\le\beta N$. Using $K=2CS^2\ge 2CS$ and \eqref{eq:choice-C}, we get
$$
a-\alpha N=(q-p)^2K-2\alpha q^2K-2\alpha S=\gamma_1K-2\alpha S\ge (2C\gamma_1-2\alpha)S\ge0,
$$
and
$$
\beta N-W=2\beta q^2K+2\beta S-(q^2-p^2)K-S=\gamma_2K-(1-2\beta)S\ge \bigl(2C\gamma_2-(1-2\beta)\bigr)S\ge0. 
$$
\end{proof}

\subsection{Correctness of the reduction}

\begin{lemma}\label{lem:forward}
If $D$ has a bipartite realization, then $D'$ has a bipartite realization.
\end{lemma}
\begin{proof}
Let $G$ be a bipartite realization of $D$ with vertex classes $A$ and $B$. The degree sum in each class is $\frac{S}{2}$, and every degree is at least $1$, hence $|A|,|B|\le\frac S2\le S$. Add $S-|A|$ isolated vertices to $A$ and $S-|B|$ isolated vertices to $B$; we obtain a bipartite graph $G_0$ with vertex classes $U$ and $V$, $|U|=|V|=S$, whose degree sequence is $(d'_1,\ldots,d'_{2S})$ (we used $|A|+|B|=m$, so exactly $2S-m$ isolated vertices are added). 

Now add four new sets of vertices: $\Lambda_U$ and $\Lambda_V$ of size $L$ each, and $\Sigma_U$ and $\Sigma_V$ of size $s$ each. The vertex classes of the new graph $G'$ will be $U\cup\Lambda_U\cup\Sigma_U$ and $V\cup\Lambda_V\cup\Sigma_V$. In addition to the edges of $G_0$, add
\begin{itemize}
    \item all edges between $\Lambda_U$ and $\Lambda_V$;
    \item all edges between $\Lambda_U$ and $V$, and all edges between $\Lambda_V$ and $U$;
    \item a semi-regular bipartite graph between $\Lambda_U$ and $\Sigma_V$ in which every vertex of $\Lambda_U$ has degree $c$ and every vertex of $\Sigma_V$ has degree $a$, and similarly between $\Lambda_V$ and $\Sigma_U$.
\end{itemize}
The semi-regular graphs exist by Lemma~\ref{lem:biregular}, since $Lc=sa$ by \eqref{eq:I2}, $c\le s$ and $a\le L$ by \eqref{eq:ineq-param}.
Clearly, $G'$ is a simple bipartite graph. A vertex of $\Lambda_U$ has degree $L+S+c=b+S=W$ by \eqref{eq:I1}; a vertex $u\in U$ has degree $d_{G_0}(u)+L$; and a vertex of $\Sigma_U$ has degree $a$. By symmetry, the same holds on the other side. Hence $G'$ is a bipartite realization of $D'$.
\end{proof}

\begin{lemma}\label{lem:backward}
If $D'$ has a bipartite realization, then $D$ has a bipartite realization.
\end{lemma}
\begin{proof}
Let $G'$ be a bipartite realization of $D'$ with vertex classes $X$ and $Y$. For $Z\in\{X,Y\}$, let $\Lambda_Z$, $O_Z$ and $\Sigma_Z$ denote the sets of large, original and small vertices in $Z$, respectively, and let
$$
o_Z:=|O_Z|,\qquad s_Z:=|\Sigma_Z|,\qquad \sigma_Z:=\sum_{v\in O_Z}\hat d(v).
$$
Then $|\Lambda_X|+|\Lambda_Y|=2L$, $o_X+o_Y=2S$, $s_X+s_Y=2s$ and $\sigma_X+\sigma_Y=S$. By symmetry, we may assume that $|\Lambda_X|\ge|\Lambda_Y|$, and we write $|\Lambda_X|=L+t$ and $|\Lambda_Y|=L-t$ with an integer $t\ge 0$.

The number of edges of $G'$ is half of the sum of the entries of $D'$, namely,
$$
E:=\frac{2LW+\sum_{v}(L+\hat d(v))+2sa}{2}=LW+SL+\frac{S}{2}+sa,
$$
where the sum is over the $2S$ original vertices. Since $G'$ is bipartite, the degree sum of $X$ and the degree sum of $Y$ are both equal to $E$.

\medskip\noindent\textbf{Step 1: $t\le\sqrt{S/2}$.} Apply Observation~\ref{obs:counting} to $A:=\Lambda_X$. For $y\in\Lambda_Y$, we use $\min\{L+t,W\}\le L+t$; for $y\in O_Y$, we have $\min\{L+t,L+\hat d(y)\}=L+\min\{t,\hat d(y)\}\le L+\hat d(y)$; and for $y\in\Sigma_Y$, we have $\min\{L+t,a\}=a$ since $a<L$. Hence
$$
(L+t)W\le(L-t)(L+t)+o_YL+\sigma_Y+s_Ya.
$$
Since the degree sum of $Y$ is $E$, we have $o_YL+\sigma_Y+s_Ya=E-(L-t)W$, and therefore
$$
(L+t)W\le L^2-t^2+E-(L-t)W,\qquad\text{that is,}\qquad t^2\le L^2+E-2LW.
$$
Using $W=b+S$ and \eqref{eq:I3},
$$
E-2LW=SL+\frac S2+sa-LW=SL+\frac S2+sa-Lb-LS=\frac S2-L^2.
$$
Hence $t^2\le\frac{S}{2}$.

\medskip\noindent\textbf{Step 2: $t=0$.} The numbers $a$, $L$ and $b$ are multiples of $K$, and hence of $2S^2$, and $W=b+S$. Comparing the degree sums of $X$ and $Y$ modulo $2S^2$, we obtain
$$
(L+t)S+\sigma_X\equiv(L-t)S+\sigma_Y\pmod{2S^2},\qquad\text{that is,}\qquad 2tS-(\sigma_Y-\sigma_X)\equiv0\pmod{2S^2}.
$$
Assume, for a contradiction, that $t\ge 1$. Since $|\sigma_Y-\sigma_X|\le S$, and $t\le\sqrt{S/2}$ by Step 1, we have
$$
0<S\le 2tS-(\sigma_Y-\sigma_X)\le 2S\sqrt{S/2}+S=S\left(\sqrt{2S}+1\right)<2S^2,
$$
where the last inequality holds since $\sqrt{2S}+1<2S$ for $S\ge2$. This contradicts the congruence above. Therefore, $t=0$, that is, $|\Lambda_X|=|\Lambda_Y|=L$.

\medskip\noindent\textbf{Step 3: the structure of $G'$.} Apply Observation~\ref{obs:counting} to $A:=\Lambda_X$ (now $|A|=L$). For $y\in\Lambda_Y\cup O_Y$, we have $\deg(y)\ge L$, and for $y\in\Sigma_Y$, we have $\deg(y)=a<L$. Hence
$$
LW\le L\cdot L+o_YL+s_Ya.
$$
Similarly, applying Observation~\ref{obs:counting} to $\Lambda_Y$ gives $LW\le L^2+o_XL+s_Xa$. Adding these two inequalities, we obtain
$$
2LW\le 2L^2+2SL+2sa.
$$
However, by $W=b+S$ and \eqref{eq:I3}, $2LW=2Lb+2LS=2L^2+2sa+2SL$; that is, equality holds in the sum, and hence in both inequalities. By the characterization of equality in Observation~\ref{obs:counting}:
\begin{itemize}
    \item every vertex of $O_Y$ is adjacent to all vertices of $\Lambda_X$;
    \item every vertex of $\Sigma_Y$ has all its neighbors in $\Lambda_X$;
\end{itemize}
and symmetrically, every vertex of $O_X$ is adjacent to all vertices of $\Lambda_Y$, and every vertex of $\Sigma_X$ has all its neighbors in $\Lambda_Y$.

Now let $x\in O_X$. Its degree is $L+\hat d(x)$, and its neighbors are in $Y=\Lambda_Y\cup O_Y\cup\Sigma_Y$. It is adjacent to all $L$ vertices of $\Lambda_Y$, and it is not adjacent to any vertex of $\Sigma_Y$, since the neighbors of the vertices of $\Sigma_Y$ are in $\Lambda_X$. Therefore, $x$ has exactly $\hat d(x)$ neighbors in $O_Y$. Symmetrically, every $y\in O_Y$ has exactly $\hat d(y)$ neighbors in $O_X$. Consequently, the subgraph of $G'$ induced by $O_X\cup O_Y$ is a bipartite graph with degree sequence $(d'_1,\ldots,d'_{2S})$. Deleting its $2S-m$ vertices with $\hat d(v)=0$ (these are isolated in the induced subgraph), we obtain a bipartite realization of $D$.
\end{proof}

\begin{proof}[Proof of Theorem~\ref{thm:hard}]
The map $D\mapsto D'$ (together with the preprocessing) is computable in polynomial time and its output is an instance of $\BDRab$ by Lemma~\ref{lem:size}. By Lemmas~\ref{lem:forward} and \ref{lem:backward}, $D'$ has a bipartite realization if and only if $D$ has one. Hence, if \BDR{} is NP-hard, then so is $\BDRab$, and since $\BDRab$ is in NP, it is NP-complete.
\end{proof}

\begin{remark}
In Construction~\ref{constr}, the original part of the sequence is padded to $S$ original vertices per vertex class, although $\frac{S}{2}$ would suffice for the forward direction. The larger padding makes the degree $W$ of the large vertices congruent to $S$ (rather than to $\frac S2$) modulo $2S^2$, and this is exactly what excludes an unbalanced distribution of the large vertices in Step~2 of the proof of Lemma~\ref{lem:backward}. With padding $\frac{S}{2}$ only, there are instances in which the large vertices are distributed unequally between the two vertex classes in some realization of $D'$ (for example, when $D$ contains many vertices of degree $1$).
\end{remark}

\section{Concluding remarks}\label{sec:concl}

We showed that $\BDRab$ is in P when $\alpha(1-\beta)\ge\beta^2$, and that it is NP-complete when $\sqrt{2\alpha}+\sqrt{1-2\beta}<1$ and $\beta<\frac12$, provided that \BDR{} itself is NP-complete. The two regions are separated by the region
$$
\frac{\sqrt{\alpha(\alpha+4)}-\alpha}{2}<\beta\le\sqrt{2\alpha}-\alpha,
$$
see Figure~\ref{fig:regions}. By Remark~\ref{rem:sharp}, in this intermediate region the existence of an equal-sum split is no longer sufficient for potential bipartiteness, so a polynomial algorithm (if it exists) must use further information. It would be interesting to close this gap, and, of course, to settle the complexity of \BDR{} itself.

\section*{Declaration of generative AI and AI-assisted technologies in the manuscript preparation process}

During the preparation of this work, the author used Claude Opus~5.5 (Anthropic), a generative AI model, as an assistant. The research problem, the main results, and the overall approach of the paper are the author's. In particular, the parameter bounds in Theorems~\ref{thm:main-structural} and~\ref{thm:hard} were found by the author, as was the approach of proving tractability via the Gale--Ryser inequalities for equal-sum splits and conditional hardness via a padding reduction. An earlier version of the manuscript (arXiv:2512.17709) contained miscalculations and gaps in the proofs. The author provided the tool with this earlier version, comments received on it, and a partially revised \LaTeX{} source, and used the tool to correct and complete the proofs and to revise the text.
 Specifically, the tool
(i) proposed the direct proof of Theorem~\ref{thm:main-structural}, based on the estimates \eqref{eq:L-bound} and \eqref{eq:R-bound} and a convexity argument, which replaces the earlier argument based on least balanced degree sequences, and suggested Remark~\ref{rem:sharp};
(ii) pointed out, with a counterexample, that the reduction in the earlier version did not force the large vertices to be distributed equally between the two vertex classes, and proposed the corrected Construction~\ref{constr}, with a larger padding of the original vertices, together with the proofs of Lemmas~\ref{lem:forward} and~\ref{lem:backward};
(iii) simplified the algebraic statements into Lemmas~\ref{lem:tractable-alg} and~\ref{lem:hard-alg};
(iv) tested the proposed arguments numerically with computer programs, namely exhaustive checks of extremal equal-sum splits for Theorem~\ref{thm:main-structural} and of the reduction on small instances; these computations are not part of any proof in the paper;
(v) corrected and completed the bibliographic data; and
(vi) drafted the \LaTeX{} source of the revised text and the \texttt{pgfplots} code of Figure~\ref{fig:regions}.
After using this tool, the author read the entire manuscript, checked all statements and proofs, edited the content as needed, and takes full responsibility for the content of the published article.

\section*{Acknowledgement} The work was supported by the Hungarian NKFIH grant K132696 and by the European Union project RRF2.3.1-21-2022-00006 within the framework of Health Safety National Laboratory Grant no RRF-2.3.1-21-2022-00006.

\end{document}